# The Dimension of One-step Idempotent Right Modular Quasigroups

R. A. R. Monzo[1]

**Abstract.** We prove that one-step, idempotent right modular groupoids are quasigroups. The dimension of such quasigroups is defined and the Cayley tables for all such quasigroups of dimensions 2, 3 and 4 are determined.

## 1. Introduction.

Redei defined finite one-step non-commutative structures as those finite structures all of whose proper substructures (of the same kind) are commutative. "The determination of such structures is a difficult task and has only been solved for groups, rings and semigroups" [3]. Since then (1967) we know of no further algebraic structures for which this problem has been solved.

This paper offers some thoughts and results that may lead to the determination of all finite one-step, idempotent right modular groupoids. In Section 2 the properties of (one-step) idempotent right modular groupoids are listed. We prove that such structures are, in fact, quasigroups and call them "one-step *W*-quasigroups" or simply "*W*-quasigroups". We prove that any two distinct elements in such a quasigroup generate it.

In Section 3, using the facts that *W*-quasigroups are non-commutative and are generated by any two distinct elements, we define the dimension of a *W*-quasigroup. In Section 4 we determine the Cayley tables of *W*-quasigroups of dimensions 2 and 3, showing that there is one *W*-quasigroup of dimension 2 (of order 4) and two of dimension 3 (of orders 5 and 11).

In Section 5 we list the Cayley tables of *W*-quasigroups of dimension 4. They are of orders 9, 11, 16, 29 and two (non-isomorphic ones) of order 19. The final Section 6 offers steps that may be able to be used to determine all one-step, idempotent right modular quasigroups, of any dimension.

Redei gave an example of an infinite one-step non-commutative semigroup [3, p. 798]. He asked whether there are infinite, one-step non-commutative groups. As far as we are aware, it is an open question as to whether one-step idempotent right modular quasigroups are finite.

## 2. Properties of one-step, idempotent right modular quasigroups

Redei defined a non-commutative structure as "… *one-step non-commutative* when its proper substructures (of the same kind) are all commutative." [3]. We modify this definition as follows:

**Definitions.** If *V* is a variety of groupoids then a **[*finite*] *one-step V-groupoid*** is a [finite] non-commutative groupoid $G \in V$ such that any proper subgroupoid of G is commutative. We define ***the variety W of idempotent, right modular groupoids*** as $W = [x = xx\,;\, xy.z = zy.x]$.

**Proposition A (properties of idempotent right modular groupoids):**
*Any groupoid $G \in W$ satisfies the following identities*:
- (1) $x = x^2$          (*idempotency*)
- (2) $xy \cdot z = zy \cdot x$     (*right modularity*)
- (3) $xy \cdot zw = xz \cdot yw$    (*mediality*)
- (4) $x \cdot yz = xy \cdot xz$     (*left distributivity*)
- (5) $xy \cdot z = xz \cdot yz$     (*right distributivity*)
- (6) $xy.x = x.yx$       (*elasticity*).

We now assume that **G is a one-step *W*-groupoid**. So henceforth G satisfies Proposition A, unless otherwise stated.

**Definition.** $\langle e, f \rangle$ denotes the subgroupoid generated by the elements *e* and *f* of the groupoid G. Also, *xyx* denotes $xy.x$ ($= x.yx$) for any *x* and *y* in G.

**Proposition 1.** *There exists* $\{e, f\} \subseteq G$ *(*$e \neq f$*) such that* $G = \langle e, f \rangle$.

Proof. Since, by definition, G is not commutative, there exists $\{e, f\} \subseteq G$ such that $ef \neq fe$. It follows that $e \neq f$. It follows also that $G = \langle e, f \rangle$, or else $\langle e, f \rangle$ is proper and, by definition, commutative, which would imply that $ef = fe$. ∎

In the results below, in the rest of this section, we assume that $G \in W$, $\{e, f\} \subseteq G$ ($e \neq f$), $ef \neq fe$ and $G = \langle e, f \rangle$.

**Proposition 2.** *If* $g = ef$ *and* $h = fe$ *then* $g \neq h, h = gf, g = he, ge = eh$ *and* $hf = fg$.

Proof. Since $ef \neq fe$, $g \neq h$. Then $h = fe = ef \cdot f = gf$,
$g = ef = fe \cdot e = he$, $ge = he \cdot e = eh$ and $hf = gf \cdot f = fg$. ∎

**Proposition 3.** *If* $g = ef$ *and* $h = fe$ *then* $gf \neq fg$ *and* $G = \langle f, g \rangle$.

Proof. By Proposition 2, $h = gf$. If $gf = fg$ then $h = gf = fg = gf \cdot f = hf$. Then
$g = ef = fe \cdot e = he = hf \cdot e = ef \cdot h = gh = gh \cdot h = hg = gf \cdot g = fg \cdot g = gf = h$, a contradiction, since $g = ef \neq fe = h$. So $gf \neq fg$ and therefore $G = \langle f, g \rangle$, or else $\langle f, g \rangle$ would be proper and, thus by hypothesis, commutative, implying $gf = fg$. ∎

**Proposition 4.** $ef \cdot fe \neq fe \cdot ef$

Proof. Assume that $ef.fe = fe.ef$. Then
$ef.fe = fe.ef = (ef.fe)(fe) = (fe.ef)(fe) = fef.efe = (ef.fe)(ef) = efe.fef$, so $ef = (e.ef)fef =$
$= (fef.ef)e = (fef.e)efe = (ef.fe)efe = (fef.efe)efe = efe.fef = fe.ef$ and similarly, $fe = ef.fe$.
Hence, $ef = fe.ef = ef.fe = fe$, a contradiction. ∎

**Corollary 5.** $G = \langle ef, fe \rangle$

**Theorem 6.** *For any $\{x, y\} \subseteq G$ there exists $\{t, s\} \subseteq G$ such that $xt = y$ and $sx = y$.*

Proof. By Corollary 5, $G = \langle ef, fe \rangle$. Let $x$ consist of the letters $x_1, x_2, ..., x_n$, ordered as they appear in $x \in G = \langle ef, fe \rangle$, with each $x_i \in \{ef, fe\}$ ($i \in 1, 2, ..., n$). We form words $x_{ef}$ and $x_{fe}$ as follows. For $j \in \{1, 2, ..., n\}$ if $x_j = ef$ define $q_j = ef$ and $k_j = f$. If $x_j = fe$ then define $q_j = e$ and $k_j = fe$.

Form the word $x_{ef}$ using the ordering $q_1, q_2, ..., q_n$ with the same bracketing as that in $x \in \langle ef, fe \rangle$, where the ordering of the letters in $x$ is $x_1, x_2, ..., x_n$. Form the word $x_{fe}$ using the ordering $k_1, k_2, ..., k_n$ with the same bracketing as that in $x \in \langle ef, fe \rangle$, where the ordering of the letters in $x$ is, again, $x_1, x_2, ..., x_n$. Then, using distributivity, it is clear that $x \cdot x_{ef} = ef$ and $x \cdot x_{fe} = fe$.

Suppose $v_1, v_2, ..., v_m$ is an ordering of the terms $ef$ and $fe$ as they appear in $y \in \langle ef, fe \rangle$. Define a word $t$ with ordered sequence $t_1, t_2, ..., t_m$ as follows. For $j \in \{1, 2, ..., m\}$, if $v_j = ef$ then $t_j = x_{ef}$ and if $v_j = fe$ then $t_j = x_{fe}$. We form the word $t$ by bracketing the $t_j$'s in the same way the $v_j$'s are bracketed in $y \in \langle ef, fe \rangle$. So, using distributivity, clearly $xt = y$ and $y = yy = xt.y = yt.x$ ∎

**Theorem 7.** $G$ *is right cancellative*

Proof. Suppose that $ax = bx$ for some $\{a, b, x\} \subseteq G$. By Theorem 6, there exists $s \in G$ such that $sx = b$. Then, $b = sx = sx.b = bx.s = ax.s = sx.a = ba = ba.a = ab$. So, $b = ab$. Also by Theorem 6, $tx = a$, for some $t \in G$. Hence, $a = tx = tx.a = ax.t = bx.t = tx.b = ab = b$. ∎

**Theorem 8.** $G$ *is left cancellative*

Proof. Suppose that $xa = xb$. Then $ax = xa.a = xb.a = ab.x$ and by Theorem 7, $a = ab = ab.b = ba$. Then, $bx = xb.b = xa.b = ba.x$ and so $b = ba = a$. ∎

**Corollary 9.** $G$ *is a quasigroup*

Proof. Using Theorems 7 and 8, the elements $t$ and $s$ of Theorem 6 are unique. Hence, $G$ is right and left solvable and, by definition, a quasigroup. ∎

**Definition.** A groupoid is ***nowhere commutative*** if $xy = yx$ implies $x = y$.

**Proposition 10.** *Idempotent, right modular, cancellative groupoids are nowhere commutative*

Proof. If $xy = yx$ then $xy = yx = xy.y = yx.x$ and, by cancellativity, $x = xy = yx = y$. ∎

**Corollary 11.** G *is nowhere commutative*

**Theorem 12.** G *is generated by any two of its distinct elements*

Proof. Since $xy = yx$ if and only if $x = y$, any two distinct elements must generate all of G. ∎

**Corollary 13.** *The only proper subgroupoids of* G *are its singleton sets.*

**Theorem 14.** G *is right and left simple*

Proof. Suppose that I is a proper right or left ideal of G. Then $I \cdot I \subseteq I$ and, by Corollary 12, I = $\{x\}$. But then, for any $y \in G$, either $xy \in I = \{x\}$ or $yx \in I = \{x\}$. By cancellation, $x = y$ and so G = I, a contradiction. ∎

**Theorem 15.** *A non-commutative, idempotent, right modular groupoid* S *is one-step if and only if it is generated by any two of its distinct elements.*

Proof. ($\Rightarrow$) Firstly, if S is one-step then, by Theorem 11, it is generated by any two of its distinct elements. ($\Leftarrow$) Conversely, If S is generated by any two of its distinct elements then any proper subgroupoid must consist of a single (idempotent) element. Hence, every proper subgroupoid of S is commutative and S is one-step. ∎

## 3. The dimension of one-step *W*-quasigroups

We have seen in Proposition 3 that, for $G \in W$, $G = \langle f, g \rangle$, where $g = ef$. This sets the stage for a definition of the "dimension" of a one-step *W*-quasigroup, as follows.

**Definition.** By a **word** $x \in G$ we mean a finite product of terms, each equal to *f* or *ef*, bracketed in any meaningful manner. The **length of a word** $x \in G = \langle f, g \rangle$ is the number of such terms in a particular designation of the word $x$ and is denoted by $|x|$, while $\|x\|$ denotes the minimum such length of the word $x$. The minimum length of the element $e$ is called **the dimension of** G, denoted by $\|e\|$. (Note: $\|e\| \neq 1$.)

## 4. One-step *W*-quasigroups of dimension two and three

**Theorem 16.** *If* $\|e\| = 2$ *then* G *is isomorphic to the following quasigroup.*

| G   | *e*  | *f*  | *ef* | *fe* |
|-----|------|------|------|------|
| *e* | e    | ef   | fe   | f    |
| *f* | fe   | f    | e    | ef   |
| *ef*| f    | fe   | ef   | e    |
| *fe*| ef   | e    | f    | fe   |

Proof. If $\|e\|=2$ then $e \in \{f, ef, fef, ef \cdot f = fe\}$. Cancellativity implies $e = fef$. Then $e.ef = fef.ef = ef.f = fe = efe.e$ and so $f = efe$. Also, $ef = fef.f = f.fe$. Then $ef.fe = fef.e = e$ and $fe.ef = efe.f = f$. It is then straightforward to prove that G is idempotent, right modular and generated by any two distinct elements. By Theorem 15, G is one-step. ■

Note that the quasigroup of Theorem 15 is $T_4$, the basic building block of idempotent, right modular, anti-rectangular groupoids [1, p. 88]. It follows from Theorem 1 that $T_4$ is the only one-step *W*-quasigroup that satisfies the identity $xy.x = y$ (anti-rectangularity, [2]).

**Lemma 17.** $\|e\|=3$ *if and only if* $e \in \{fef \cdot f, f \cdot fef\}$

Proof. ($\Leftarrow$) If $e \in \{fef \cdot f, f \cdot fef\}$ then $\|e\|=3$, since from Theorem 15, $\|e\|=2$ implies $e = fef$ and, by cancellation, $f = fef$. But $f = fef$ implies $e = f$, a contradiction.

($\Rightarrow$) $\|e\|=3$ implies $e \in \{f, g, fg, fgf, gf, gfg, f \cdot fg, g \cdot gf\} =$
$= \{f, ef, fef, fef \cdot f, fe, fe \cdot ef, f \cdot fef, ef \cdot fe\}$. But since $e \notin \{f, ef, fe\}$ and $e = ef.fe$ implies $e = fef.e$ implies $e = fef$ implies $\|e\|=2$, we can assume that $e \in \{fef \cdot f, fe \cdot ef, f \cdot fef\}$. We show that $e = fef.f$ if and only if $e = fe.ef$. If $e = fef.f$ then $e = f.fe = (f.fe)e = fe.ef$. Conversely, $e = fe.ef$ implies $e = (fe.e)fef = ef.fef = (fef.f)e$, which by cancellation implies $e = fef.f$. Therefore, $e \in \{fef \cdot f, f \cdot fef\}$. ■

Note that the proof of Theorem 18 below relies heavily on Proposition A and Theorems 7 and 8, as do the omitted calculations that determine the Cayley Tables 3 – 10 below.

**Theorem 18.** $\|e\|=3$ *implies* G *is one of the following two quasigroups*:
*If* $e = fef \cdot f$ *then* G *has the following Cayley table*:

**Table 1.** $e = fef \cdot f$

| G    | *e*  | *f*  | *ef* | *fe* | *fef* |
|------|------|------|------|------|-------|
| *e*  | e    | ef   | f    | fef  | fe    |
| *f*  | fe   | f    | fef  | e    | ef    |
| *ef* | fef  | fe   | ef   | f    | e     |
| *fe* | ef   | fef  | e    | fe   | f     |
| *fef*| f    | e    | fe   | ef   | fef   |

If $e = f.fef$ then G has the following Cayley table, where the ordered pair
(1,2,3,4,5,6,7,8,9,10,11) = (*e, ef, fe, fef, f.fe, e.ef, efe, ef.fe, fe.ef, e(f.fe), e.efe*).

**Table 2.** $e = f \cdot fef$

| G  | 1  | 2  | 3  | 4  | 5  | 6  | 7  | 8  | 9  | 10 | 11 |
|----|----|----|----|----|----|----|----|----|----|----|----|
| 1  | 1  | 6  | 7  | 5  | 10 | 3  | 11 | 9  | 4  | 8  | 2  |
| 2  | 7  | 2  | 8  | 9  | 1  | 5  | 10 | 4  | 11 | 6  | 3  |
| 3  | 2  | 9  | 3  | 10 | 6  | 11 | 8  | 1  | 7  | 5  | 4  |
| 4  | 8  | 3  | 6  | 4  | 7  | 9  | 1  | 11 | 10 | 2  | 5  |
| 5  | 9  | 10 | 4  | 8  | 5  | 7  | 3  | 6  | 2  | 11 | 1  |
| 6  | 11 | 7  | 10 | 2  | 4  | 6  | 9  | 3  | 5  | 1  | 8  |
| 7  | 6  | 5  | 2  | 11 | 8  | 4  | 7  | 10 | 1  | 3  | 9  |
| 8  | 5  | 11 | 9  | 7  | 3  | 1  | 2  | 8  | 6  | 4  | 10 |
| 9  | 10 | 8  | 1  | 3  | 11 | 2  | 4  | 5  | 9  | 7  | 6  |
| 10 | 4  | 1  | 11 | 6  | 9  | 8  | 5  | 2  | 3  | 10 | 7  |
| 11 | 3  | 4  | 5  | 1  | 2  | 10 | 6  | 7  | 8  | 9  | 11 |

Proof. We know from Lemma 17 that $\|e\| = 3$ *implies* $e \in \{fef \cdot f, f \cdot fef\}$ and, from the proof of Lemma 16, that *e = fef.f* if and only if *e = fe.ef*.

Assume *e = fef.f*, which implies *e = fe.ef*. So, *ef.fe = (fe.ef)(f.fe) = (fef.fe) fe = fe.fef = f(e.ef)*. Since *e.ef = (fe.ef)ef = ef.fe = f(e.ef)*, we have that *f = e.ef* and *ef.fe = f(e.ef) = f*.

Also, *fef = f.ef = (e.ef)ef = efe*, and therefore $fef \neq e$. So $fef \notin \{e, f, ef, fe\}$.
Then, *f.fe = (ef.fe)fe = fe.ef = e*. It is then straightforward to calculate that **the Cayley table for G is as indicated in Table 1 of Theorem 18 when** *e = fef.f*

Now suppose that *e = f.fef*. Note that, since
*(e.efe)e =[(e.ef)e]e = e(e.ef) = (f.fef)(e.ef) =(fe)(fef.ef) = (fe)(ef.f) = fe*. We have *f = e.efe*
Therefore, the second row of Table two has the following products:

| G | e | ef    | fe  | fef | f(fe)  | e(ef) | efe   | (ef)(fe) | (fe)(ef) | e(f.fe) | e.efe |
|---|---|-------|-----|-----|--------|-------|-------|----------|----------|---------|-------|
| e | e | e(ef) | efe |     | e(f.fe)|       | e.efe |          |          |         | ef    |

Then, *e.fef = (f.fef)(fef) = fef.f = f.fe* , *e(e.ef) = [(e.ef)e]e = (e.efe)e = fe* ,
*e(ef.fe) = e[(ef.f)(efe)] = e(fe.efe) = (efe)(e.efe) = efe.f = fe.ef* ,
*e(fe.ef) = (f.fef)(fe.ef) = (f.fe)(fef.ef) = (f.fe)(fe) = fef*
and *e[e(f.fe)] = (f.fef) [e(f.fe)] = (fe)[f(ef.fe)] = f[e(ef.fe)] = f (fe.ef) = (e.efe)(fe.ef) =*
*= (efe)[(efe)(ef)] = (efe)[(e.efe)(fe)] = (efe)(f.fe) = ef.fe*
So the products in the ***second*** row of Table two are valid.
The following products in the third row of Table two are, clearly, as follows:

| G  | e   | ef  | fe     | fef | f(fe) | e(ef) | efe    | (ef)(fe) | (fe)(ef) | e.(f.fe) | e.efe |
|----|-----|-----|--------|-----|-------|-------|--------|----------|----------|----------|-------|
| ef | efe | ef  | (ef)(fe) |     |       |       | e(f.fe) |          |          |          | fe    |

Then, *(ef)(fef) = [(fef)(ef)](ef) = (ef.f)(ef) = fe.ef* ,
*(ef)(f.fe) = (ef)(fef.f) = (e.fef)f = (f.fe)f = f.fef = e* ,
*(ef)(e.ef) = e.fef = f.fe* ,
*(ef)(ef.fe) = (ef)(fef.e) = (e.fef)(fe) = (f.fe)(fe) = fef*

$(ef)(fe.ef) = (fe.e)(fe.ef) = (fe)(e.ef) = [(e.ef)e]f = (e.efe)f = e.efe$ and
$(ef)[e(f.fe)] = \{[e(f.fe)]f\}e = \{(ef)[(f.fe)f]\}e = [(ef)(f.fef)]e = efe.e = e.ef$

So the products in the *third* row of Table two are valid
The following products in the fourth row of Table two are, clearly, as follows:

| G  | e  | ef      | fe | fef | f(fe) | e(ef) | efe | (ef)(fe) | (fe)(ef) | e(f.fe) | e.efe |
|----|----|---------|----|-----|-------|-------|-----|----------|----------|---------|-------|
| fe | ef | (fe)(ef)| fe |     |       |       |     |          |          |         | fef   |

Then,
$(fe)(fef) = f(e.ef) = (e.efe)(e.ef) = [(e.ef)(efe)]e = [e(ef.fe)]e = (fe.ef)e = ef.efe = e(f.fe)$,
$(fe)(f.fe) = f.efe = (e.efe)(efe) = efe.e = e.ef$,
$(fe)(e.ef) = [(e.ef)e]f = (e.efe)f = e.efe$,
$(fe)(efe) = (fe.e)(fe.fe) = ef.fe$
$(fe)(ef.fe) = (fe)(fef.e) = (f.fef)e = e$
$(fe)(fe \cdot ef) = (fe)[(ef.fe)(fe)] = \{(fe)[ef \cdot fe]\} \cdot (fe) = efe$ and
$(fe)[e(f.fe)] = (fe)(ef.efe) = (fe.ef)(fe.efe) = (fe.ef)(ef.fe) = [(fe)(ef.fe)][(ef)(fe.ef)] = e.fef = f.fe$

So the products in the *fourth* row of Table two are valid
The following products in the fifth line of Table two are clearly valid:

| G   | e     | ef | fe | fef | f(fe) | e(ef)  | efe | (ef)(fe) | (fe)(ef) | e(f.fe) | e.efe |
|-----|-------|----|----|-----|-------|--------|-----|----------|----------|---------|-------|
| fef | ef·fe | fe |    | fef |       | fe·ef  |     |          |          |         | f·fe  |

Then, $fef.fe = (fe)(f.fe) = e.ef$,
$(fef)(f.fe) = f(ef.fe) = (e.efe)(ef.fe) = (e.ef)(efe.fe) = (e.ef)(ef) = efe$
$fef.efe = (f.efe)(ef.efe) = (e.ef)[e(f.fe)] = \{[e(f.fe)](ef)\}e = [e(f.fef)]e = e$
$fef \cdot (ef \cdot fe) = fef \cdot (fef \cdot e) = [(fef \cdot e)f](fe) = (fe \cdot fef)(fe) = [(f \cdot fef)(e \cdot fef)](fe) =$
$[e(f \cdot fe)](fe) = (ef)[(f \cdot fe)e] = (ef)(fe \cdot ef) = e.efe$, $fef \cdot (fe \cdot ef) = (f \cdot fe)(ef) = fe \cdot fef = e(f \cdot fe)$ and
$fef \cdot [e \cdot (f \cdot fe)] = (fe \cdot e)[f \cdot (f \cdot fe)] = ef \cdot ef = ef$.
So the products in the *fifth* row of Table two are valid.
The following products in the sixth row of Table two are clearly valid:

| G     | e     | ef     | fe  | fef     | f(fe) | e(ef) | efe | (ef)(fe) | (fe)(ef) | e(f.fe) | e.efe |
|-------|-------|--------|-----|---------|-------|-------|-----|----------|----------|---------|-------|
| f(fe) | fe.ef | e(f.fe)| fef | (ef)(fe)| f(fe) | efe   | fe  | e(ef)    | ef       |         | e     |

Then, $f(fe) \cdot e[f(fe)] = fe \cdot [(fe)(f \cdot fe)] = fe \cdot (e \cdot ef) = [e \cdot efe]f = f$

So the products in the *sixth* row of Table two are valid.
The following products in the seventh row of Table two are clearly valid:

| G    | e      | ef  | fe     | fef | f(fe) | e(ef) | efe   | (ef)(fe) | (fe)(ef) | e(f.fe) | e.efe   |
|------|--------|-----|--------|-----|-------|-------|-------|----------|----------|---------|---------|
| e.ef | e.efe  | efe | e(f.fe)| ef  | fef   | e(ef) | (fe)(ef)| fe     | f(fe)    |         | (ef)(fe)|

Then, $(e \cdot ef) \cdot e[f(fe)] = \{e[f(fe)] \cdot ef\}e = e$.

So the products in the *seventh* row of Table two are valid.
The following products in the eighth row of Table two are clearly valid:

| G   | e     | ef    | fe | fef    | f(fe)  | e(ef) | efe | (ef)(fe) | (fe)(ef) | e(f.fe) | e(efe) |
|-----|-------|-------|----|--------|--------|-------|-----|----------|----------|---------|--------|
| efe | e(ef) | f(fe) | ef | e(efe) | (ef)(fe)| fef  | efe |          |          |         |        |

Then, $efe \cdot (ef \cdot fe) = (ef \cdot fe)(fe) \cdot e = (fe \cdot ef) \cdot e = [e(ef)](fe) = e[f(fe)]$,
$efe \cdot (fe \cdot ef) = (efe \cdot fe)(efe \cdot ef) = (ef)(f \cdot fe) = e$,
$efe \cdot e[f(fe)] = e[fe \cdot f(fe)] = e \cdot (e \cdot ef) = [(e \cdot ef)e]e = fe$ and
$efe \cdot e(efe) = efe \cdot f = fe \cdot ef$.

So the products in the *eighth* row of Table two are valid:
The following products in the ninth row of Table two are clearly valid:

| G | e | ef | fe | fef | f(fe) | e(ef) | efe | (ef)(fe) | (fe)(ef) | e(f.fe) | e(efe) |
|---|---|----|----|-----|-------|-------|-----|----------|----------|---------|--------|
| ef.fe | f(fe) | e(efe) | fe.ef |  | fe |  |  | (ef)(fe) | e(ef) |  | e[f(fe)] |

Then, $(ef \cdot fe)(fef) = (fef \cdot fe)(ef) = (e \cdot ef)(ef) = efe$,
$(ef \cdot fe)(e \cdot ef) = efe \cdot (fe \cdot ef) = e$,
$(ef \cdot fe)(efe) = [(efe)(fe)](ef) = ef$ and
$(ef \cdot fe) \cdot e[f(fe)] = (efe)(f \cdot efe) = (efe)(e \cdot ef) = fef$.

So the products in the *ninth* row of Table two are valid:
The following products in the tenth row of Table two are clearly valid:

| G | e | ef | fe | fef | f(fe) | e(ef) | efe | (ef)(fe) | (fe)(ef) | e(f.fe) | e(efe) |
|---|---|----|----|-----|-------|-------|-----|----------|----------|---------|--------|
| fe.ef | e[f(fe)] | (ef)(fe) | e | fe |  | ef |  | f(fe) | (fe)(ef) |  |  |

Then, $(fe \cdot ef) \cdot f(fe) = (fef)(ef \cdot fe) = e \cdot efe$,
$(fe \cdot ef) \cdot (efe) = [(efe)(ef)](fe) = f(fe) \cdot fe = fef$,
$(fe \cdot ef) \cdot e[f(fe)] = (ef)[(ef)(f \cdot fe)] = efe$ and
$(fe \cdot ef) \cdot [e \cdot efe] = (fe \cdot ef) \cdot f = fef \cdot fe = e \cdot ef$.

So the products in the *tenth* row of Table two are valid:
The following products in the eleventh row of Table two are clearly valid:

| G | e | ef | fe | fef | f(fe) | e(ef) | efe | (ef)(fe) | (fe)(ef) | e(f.fe) | e.efe |
|---|---|----|----|-----|-------|-------|-----|----------|----------|---------|-------|
| e(f.fe) |  |  |  |  | (fe)(ef) |  | f(fe) |  |  | e(f.fe) | efe |

Then, $e[f(fe)] \cdot e = e \cdot (fe)(ef) = fef$,
$e[f(fe)] \cdot ef = e \cdot (f \cdot fef) = e$,
$e[f(fe)] \cdot fe = (ef)[(fe)(ef)] = e(efe)$,
$e[f(fe)] \cdot fef = \{e[f(fe)] \cdot fe\} \{e[f(fe)] \cdot [e(efe)]\} = [e(efe)](efe) = e(ef)$,
$e[f(fe)] \cdot e(ef) = e(fe \cdot fef) = (efe)(e \cdot fef) = (efe)[f(fe)] = (ef)(fe)$,
$e[f(fe)] \cdot (ef)(fe) = [e(ef)](fef) = ef$ and
$e[f(fe)] \cdot (fe)(ef) = \{e[f(fe)] \cdot fe\} \{e[f(fe)] \cdot ef\} = [e(efe)]e = fe$.

So the products in the *eleventh* row of Table two are valid:
The following products in the twelfth row of Table two are clearly valid:

| G | e | ef | fe | fef | f(fe) | e(ef) | efe | (ef)(fe) | (fe)(ef) | e(f.fe) | e.efe |
|---|---|----|----|-----|-------|-------|-----|----------|----------|---------|-------|
| e(efe) | fe | fef | f(fe) | e |  |  | e(ef) |  |  |  | e.efe |

Then, $e(efe) \cdot f(fe) = (ef)(efe \cdot fe) = ef$,
$e(efe) \cdot e(ef) = e \cdot (efe \cdot ef) = e \cdot f(fe)$,
$e(efe) \cdot (ef)(fe) = (e \cdot ef)(efe \cdot fe) = (e \cdot ef)(ef) = efe$,

$e(efe) \cdot (fe)(ef) = \{e(efe) \cdot (fe)\} \{e(efe) \cdot (ef)\} = f(fe) \cdot (fef) = (ef)(fe)$ and
$e(efe) \cdot e[f(fe)] = e[(ef)(fe)] = (fe)(ef)$.

So the products in the *twelfth* row of Table two are valid. Hence, we have proved that Table two is valid.

Note that if Table 2 is valid then the elements listed in it comprise all of G, or else the elements listed must commute, a contradiction since $ef \neq fe$. Also, using Proposition A and Table 2, it is straightforward to prove that G consists of eleven distinct elements. We omit the details of that part of the proof. The calculations showing that the groupoids determined by the Cayley Tables 1 and 2 are right modular are laborious and are also omitted. By Theorem 15, we then need only prove that any two distinct elements in Table 1 or in Table 2 generate G. Since $G = \langle e, f \rangle$, we need only show that $e$ and $f$ are equal to words whose factors consist of any two distinct elements. This we proceed to do.

In Table 1 we have $f = e.ef = e.(fe.e) = e[(e.fef)e] = ef.fe = (ef)(fef.ef) = (fe)(fef)$ and $e = (ef)(f.ef) = (fe.f)f = (fef)f$. So $e$ and $f$ are equal to words whose factors consist of any two distinct elements of G in Table 1.

In Table 2, note that $ef = 2 = 1.11 = e(e.efe)$. Therefore, $f = e.efe$. Using this fact it is straightforward to prove that $e$ and $f$ are equal to words whose factors consist of any two distinct elements of G in Table 2. We omit these detailed calculations. This completes the proof of Theorem 18. ∎

## 5. One-step *W*-quasigroups of dimension four

The first step in determining one-step *W*-quasigroups when $\|e\| = 4$ is to determine those possible values of $e$, when $\|e\| = 4$. We have calculated those values by hand, when n = 4. We list them next, without proof, and along with their equivalent expressions in the generators *ef* and *fe*: (1) $e = f[f(fef)] = f[(fef.f)f] = f[(f.fe)f] = \{f[(f.fe)f]\}e = (fe)(fe.ef.fe)$
      (2) $e = f[f(ef \cdot f)] = f(f.fe) = (fe)(fe.ef)$
      (3) $e = (f \cdot fef)(ef) = (fef)(fe) = [(fef)(fe)]e = ef.fe.ef$
      (4) $e = (fef \cdot f)(ef) = fe.fef = (fe.fef)e = (ef)(ef.fe)$
      (5) $e = f[(ef)(fef)] = f(fe.ef) = [f(fe.ef)]e = (fe)(ef.efe) = (fe.ef)(fe.efe) = (fe.ef)(fef.e) =$
          $= (fe.ef)(ef.fe)$
      (6) $e = f[(ef)(ef.f)] = f(ef.fe) = [f(ef.fe)]e = (fe)(efe.ef) = (fe)(e.fef) = (fe.e)(fe.fef) =$
          $= (ef)(fe.fef) = (ef.fe)(ef.fef) = (ef.fe)[(ef.f)(ef.ef)] = (ef.fe)(fe.ef)$

Note that there are 5 ways to multiply four elements in any groupoid: $a(bc.d)$, $a(b.cd)$, $(a.bc)d$, $ab.cd$ and $(ab.c)d$. Therefore, there are 80 different ways to multiply the generators $f$ and $ef$ when $\|e\| = 4$. All products of $f$ and $ef$ of length 4, other than those listed in (1) to (6) above, are reducible to a product of length 3 or simply do not exist. For example, $e = fe.ef.fe$ $=(f.fe.f)e$ implies $e = f.fef$, which contradicts $\|e\| = 4$. The value $e = (ef)(ef.fe.ef)$ does not arise as a product of four terms, each term being from the set $\{f, ef\}$.

Without proof, we note that G consists of the following 29 elements when *e* = *f* [*f(fef)*];
namely, if the following ordered pairs are defined as equal;
(1,2,3,4,5,6,7,8,9,10,11,12,13,14,15,16,17,18,19,20,21,22,23,24,25,26,27,28,29) =
= (*e*, *f*, *ef*, *fe*, *fef*, *f.fe*, *e.ef*, *efe*, *e(f.fe)*, *ef.fe*, *fe.ef*, *f.efe*, *e.efe*, *e.fef*, *fef.efe*, *efe.fef*, *e(e.ef)*, *f(e.ef)*, *f(f.fe)*, *f(fef)*, *e(ef.fe)*, *f(ef.fe)*, *f(fe.ef)*, *(fe)(fe.ef)*, *(ef)(ef.fe)*, *f(e.fef)*, *(ef.fe)(fe.ef)*, *(fe.ef)(ef.fe)*, *f [e(e.ef)]*) . Then the Tables 3,4 and 5 below are the Cayley table of G when
(1) *e* = *f[f.fef]*. We omit the detailed but straightforward calculations.

### Table 3.  *e* = *f* (*f.fef*)

| G  | 1  | 2  | 3  | 4  | 5  | 6  | 7  | 8  | 9  | 10 |
|----|----|----|----|----|----|----|----|----|----|----|
| 1  | 1  | 3  | 7  | 8  | 14 | 9  | 17 | 13 | 27 | 21 |
| 2  | 4  | 2  | 5  | 6  | 20 | 19 | 18 | 12 | 7  | 22 |
| 3  | 8  | 4  | 3  | 10 | 11 | 15 | 14 | 9  | 26 | 25 |
| 4  | 3  | 5  | 11 | 4  | 18 | 12 | 16 | 10 | 28 | 15 |
| 5  | 10 | 6  | 4  | 12 | 5  | 22 | 11 | 15 | 1  | 17 |
| 6  | 11 | 20 | 18 | 5  | 23 | 6  | 24 | 4  | 17 | 12 |
| 7  | 13 | 10 | 8  | 9  | 3  | 25 | 7  | 21 | 24 | 29 |
| 8  | 7  | 11 | 14 | 3  | 16 | 10 | 22 | 8  | 29 | 9  |
| 9  | 22 | 24 | 19 | 16 | 27 | 11 | 6  | 14 | 9  | 3  |
| 10 | 14 | 18 | 16 | 11 | 24 | 4  | 19 | 3  | 25 | 10 |
| 11 | 9  | 12 | 10 | 15 | 4  | 17 | 3  | 25 | 23 | 28 |
| 12 | 16 | 23 | 24 | 18 | 26 | 5  | 27 | 11 | 15 | 4  |
| 13 | 17 | 16 | 22 | 14 | 19 | 3  | 12 | 7  | 21 | 8  |
| 14 | 21 | 15 | 9  | 25 | 10 | 28 | 8  | 29 | 18 | 26 |
| 15 | 19 | 26 | 27 | 24 | 29 | 18 | 2  | 16 | 10 | 11 |
| 16 | 29 | 17 | 25 | 28 | 15 | 1  | 9  | 26 | 5  | 23 |
| 17 | 2  | 25 | 21 | 29 | 9  | 26 | 13 | 27 | 11 | 24 |
| 18 | 25 | 22 | 15 | 17 | 12 | 7  | 10 | 28 | 20 | 1  |
| 19 | 24 | 1  | 26 | 23 | 28 | 20 | 29 | 18 | 12 | 5  |
| 20 | 15 | 19 | 12 | 22 | 6  | 14 | 4  | 17 | 13 | 7  |
| 21 | 12 | 27 | 6  | 19 | 2  | 16 | 5  | 22 | 8  | 14 |
| 22 | 27 | 28 | 29 | 26 | 25 | 23 | 21 | 24 | 4  | 18 |
| 23 | 28 | 14 | 17 | 7  | 22 | 8  | 15 | 1  | 2  | 13 |
| 24 | 26 | 7  | 28 | 1  | 17 | 13 | 25 | 23 | 6  | 20 |
| 25 | 6  | 29 | 2  | 27 | 21 | 24 | 20 | 19 | 3  | 16 |
| 26 | 23 | 8  | 1  | 13 | 7  | 21 | 28 | 20 | 19 | 2  |
| 27 | 18 | 13 | 23 | 20 | 1  | 2  | 26 | 5  | 22 | 6  |
| 28 | 20 | 9  | 13 | 21 | 8  | 29 | 1  | 2  | 16 | 27 |
| 29 | 5  | 21 | 20 | 2  | 13 | 27 | 23 | 6  | 14 | 19 |

**Table 4.** $e = f(f.fef)$.

| G | 11 | 12 | 13 | 14 | 15 | 16 | 17 | 18 | 19 | 20 |
|---|---|---|---|---|---|---|---|---|---|---|
| 1 | 22 | 29 | 2 | 12 | 24 | 5 | 4 | 6 | 26 | 19 |
| 2 | 23 | 14 | 17 | 26 | 8 | 25 | 29 | 28 | 3 | 1 |
| 3 | 16 | 28 | 29 | 19 | 23 | 2 | 6 | 27 | 1 | 24 |
| 4 | 24 | 17 | 25 | 27 | 1 | 21 | 2 | 29 | 7 | 26 |
| 5 | 18 | 7 | 28 | 24 | 13 | 29 | 27 | 26 | 8 | 23 |
| 6 | 26 | 22 | 15 | 29 | 7 | 9 | 21 | 25 | 14 | 28 |
| 7 | 14 | 26 | 27 | 22 | 18 | 6 | 12 | 19 | 23 | 16 |
| 8 | 19 | 25 | 21 | 6 | 26 | 20 | 5 | 2 | 28 | 27 |
| 9 | 2 | 10 | 8 | 20 | 25 | 1 | 23 | 13 | 15 | 21 |
| 10 | 27 | 15 | 9 | 2 | 28 | 13 | 20 | 21 | 17 | 29 |
| 11 | 11 | 1 | 26 | 16 | 20 | 27 | 19 | 24 | 13 | 18 |
| 12 | 29 | 12 | 10 | 21 | 17 | 8 | 13 | 9 | 22 | 25 |
| 13 | 6 | 9 | 13 | 5 | 29 | 23 | 18 | 20 | 25 | 2 |
| 14 | 3 | 23 | 24 | 14 | 5 | 19 | 22 | 16 | 20 | 11 |
| 15 | 21 | 4 | 3 | 13 | 15 | 7 | 1 | 8 | 12 | 9 |
| 16 | 10 | 20 | 18 | 3 | 6 | 16 | 14 | 11 | 2 | 4 |
| 17 | 8 | 18 | 16 | 7 | 4 | 22 | 17 | 14 | 5 | 3 |
| 18 | 4 | 13 | 23 | 11 | 2 | 24 | 16 | 18 | 21 | 5 |
| 19 | 25 | 6 | 4 | 9 | 22 | 3 | 8 | 10 | 19 | 15 |
| 20 | 5 | 8 | 1 | 18 | 21 | 26 | 24 | 23 | 9 | 20 |
| 21 | 20 | 3 | 7 | 23 | 9 | 28 | 26 | 1 | 10 | 13 |
| 22 | 9 | 5 | 11 | 8 | 12 | 14 | 7 | 3 | 6 | 10 |
| 23 | 12 | 21 | 20 | 4 | 27 | 18 | 11 | 5 | 29 | 6 |
| 24 | 15 | 2 | 5 | 10 | 19 | 11 | 3 | 4 | 27 | 12 |
| 25 | 13 | 11 | 14 | 1 | 10 | 17 | 28 | 7 | 4 | 8 |
| 26 | 17 | 27 | 6 | 15 | 16 | 4 | 10 | 12 | 24 | 22 |
| 27 | 28 | 19 | 12 | 25 | 14 | 10 | 9 | 15 | 16 | 17 |
| 28 | 7 | 24 | 19 | 17 | 11 | 12 | 15 | 22 | 18 | 14 |
| 29 | 1 | 16 | 22 | 28 | 3 | 15 | 25 | 17 | 11 | 7 |

**Table 5.** $e = f(f.fef)$

| G | 21 | 22 | 23 | 24 | 25 | 26 | 27 | 28 | 29 |
|---|----|----|----|----|----|----|----|----|----|
| 1 | 16 | 18 | 20 | 23 | 11 | 28 | 25 | 15 | 10 |
| 2 | 13 | 9 | 15 | 10 | 21 | 11 | 16 | 24 | 27 |
| 3 | 18 | 20 | 21 | 13 | 5 | 7 | 17 | 22 | 12 |
| 4 | 23 | 13 | 9 | 8 | 20 | 14 | 22 | 19 | 6 |
| 5 | 20 | 21 | 25 | 9 | 2 | 3 | 14 | 16 | 19 |
| 6 | 1 | 8 | 10 | 3 | 13 | 16 | 19 | 27 | 2 |
| 7 | 11 | 5 | 2 | 20 | 4 | 1 | 28 | 17 | 15 |
| 8 | 24 | 23 | 13 | 1 | 18 | 17 | 15 | 12 | 4 |
| 9 | 29 | 28 | 7 | 17 | 26 | 12 | 4 | 5 | 18 |
| 10 | 26 | 1 | 8 | 7 | 23 | 22 | 12 | 6 | 5 |
| 11 | 5 | 2 | 29 | 21 | 6 | 8 | 7 | 14 | 22 |
| 12 | 28 | 7 | 3 | 14 | 1 | 19 | 6 | 2 | 20 |
| 13 | 27 | 26 | 1 | 28 | 24 | 15 | 10 | 4 | 11 |
| 14 | 4 | 6 | 27 | 2 | 12 | 13 | 1 | 7 | 17 |
| 15 | 25 | 17 | 14 | 22 | 28 | 6 | 5 | 20 | 23 |
| 16 | 12 | 19 | 24 | 27 | 22 | 21 | 13 | 8 | 7 |
| 17 | 10 | 12 | 19 | 6 | 15 | 20 | 23 | 1 | 28 |
| 18 | 6 | 27 | 26 | 29 | 19 | 9 | 8 | 3 | 14 |
| 19 | 17 | 14 | 11 | 16 | 7 | 27 | 2 | 21 | 13 |
| 20 | 2 | 29 | 28 | 25 | 27 | 10 | 3 | 11 | 16 |
| 21 | 21 | 25 | 17 | 15 | 29 | 4 | 11 | 18 | 24 |
| 22 | 15 | 22 | 16 | 19 | 17 | 2 | 20 | 13 | 1 |
| 23 | 19 | 24 | 23 | 26 | 16 | 25 | 9 | 10 | 3 |
| 24 | 22 | 16 | 18 | 24 | 14 | 29 | 21 | 9 | 8 |
| 25 | 9 | 15 | 22 | 12 | 25 | 5 | 18 | 23 | 26 |
| 26 | 14 | 11 | 5 | 18 | 3 | 26 | 29 | 25 | 9 |
| 27 | 7 | 3 | 4 | 11 | 8 | 24 | 27 | 29 | 21 |
| 28 | 3 | 4 | 6 | 5 | 10 | 23 | 26 | 28 | 25 |
| 29 | 8 | 10 | 12 | 4 | 9 | 18 | 24 | 26 | 29 |

We proceed, without the detailed proof, to show the Cayley table when (2) $e = f[f(ef \cdot f)]$. In this case there are 16 elements. Again we number them, to facilitate showing the Cayley Table, by defining the following ordered pairs to be equal:
(1,2,3,4,5,6,7,8,9,10,11,12,13,14,15,16) =
= (*e, ef, fe, fef, f.fe, efe, e.ef, ef.fe ,fe.ef, f.efe, e.efe, e(f.fe), e.fef, efe.fef, fef.efe, e(e.ef)*).

**Table 6.** $e = f[f(ef \cdot f)]$

| G  | 1  | 2  | 3  | 4  | 5  | 6  | 7  | 8  | 9  | 10 | 11 | 12 | 13 | 14 | 15 | 16 |
|----|----|----|----|----|----|----|----|----|----|----|----|----|----|----|----|----|
| 1  | 1  | 7  | 6  | 13 | 12 | 11 | 16 | 4  | 15 | 5  | 3  | 10 | 8  | 9  | 14 | 2  |
| 2  | 6  | 2  | 8  | 9  | 15 | 12 | 13 | 16 | 14 | 7  | 5  | 1  | 10 | 4  | 11 | 3  |
| 3  | 2  | 9  | 3  | 11 | 10 | 8  | 14 | 15 | 1  | 13 | 16 | 7  | 5  | 12 | 6  | 4  |
| 4  | 8  | 3  | 10 | 4  | 14 | 15 | 9  | 13 | 11 | 2  | 7  | 6  | 1  | 16 | 12 | 5  |
| 5  | 9  | 11 | 4  | 6  | 5  | 3  | 1  | 10 | 7  | 14 | 15 | 13 | 16 | 8  | 2  | 12 |
| 6  | 7  | 13 | 2  | 14 | 8  | 6  | 15 | 12 | 10 | 16 | 4  | 5  | 3  | 11 | 1  | 9  |
| 7  | 11 | 6  | 12 | 2  | 16 | 4  | 7  | 5  | 13 | 1  | 10 | 14 | 15 | 3  | 9  | 8  |
| 8  | 13 | 14 | 9  | 1  | 3  | 2  | 10 | 8  | 5  | 15 | 12 | 16 | 4  | 6  | 7  | 11 |
| 9  | 12 | 8  | 15 | 3  | 13 | 16 | 2  | 7  | 9  | 6  | 1  | 11 | 14 | 5  | 4  | 10 |
| 10 | 14 | 1  | 11 | 7  | 4  | 9  | 5  | 3  | 16 | 10 | 8  | 15 | 12 | 2  | 13 | 6  |
| 11 | 16 | 15 | 13 | 10 | 2  | 7  | 8  | 6  | 3  | 12 | 11 | 4  | 9  | 1  | 5  | 14 |
| 12 | 15 | 10 | 14 | 5  | 9  | 13 | 3  | 2  | 4  | 8  | 6  | 12 | 11 | 7  | 16 | 1  |
| 13 | 4  | 12 | 16 | 8  | 7  | 5  | 6  | 1  | 2  | 11 | 14 | 9  | 13 | 10 | 3  | 15 |
| 14 | 5  | 16 | 7  | 15 | 6  | 1  | 12 | 11 | 8  | 4  | 9  | 3  | 2  | 14 | 10 | 13 |
| 15 | 10 | 5  | 1  | 16 | 11 | 14 | 4  | 9  | 12 | 3  | 2  | 8  | 6  | 13 | 15 | 7  |
| 16 | 3  | 4  | 5  | 12 | 1  | 10 | 11 | 14 | 6  | 9  | 13 | 2  | 7  | 15 | 8  | 16 |

We now show, without proof, the Cayley Tables when (3) $e = fef.fe$, Table 7, and (4) $e = fe.fef$, Table 8.

**Table 7.** $e = fef.fe$

|        | *e*    | *f*    | *ef*   | *fe*   | *e.ef* | *efe*  | *fef*  | *ef.fe*| *fe.ef*|
|--------|--------|--------|--------|--------|--------|--------|--------|--------|--------|
| *e*    | e      | ef     | e.ef   | efe    | fef    | fe.ef  | f      | fe     | ef.fe  |
| *f*    | fe     | f      | fef    | e.ef   | efe    | e      | ef.fe  | fe.ef  | ef     |
| *ef*   | efe    | fe     | ef     | ef.fe  | f      | fef    | fe.ef  | e.ef   | e      |
| *fe*   | ef     | fef    | fe.ef  | fe     | e      | ef.fe  | efe    | f      | e.ef   |
| *e.ef* | fe.ef  | ef.fe  | efe    | fef    | e.ef   | fe     | ef     | e      | f      |
| *efe*  | e.ef   | fe.ef  | f      | ef     | ef.fe  | efe    | e      | fef    | fe     |
| *fef*  | ef.fe  | e.ef   | fe     | e      | fe.ef  | f      | fef    | ef     | efe    |
| *ef.fe*| f      | efe    | e      | fe.ef  | fe     | ef     | e.ef   | ef.fe  | fef    |
| *fe.ef*| fef    | e      | ef.fe  | f      | ef     | e.ef   | fe     | efe    | fe.ef  |

**Table 8.** $e = fe.fef$

|   | e | f | ef | fe | fef | efe | f.fe | ef.fe | fe.ef | e.ef | e.fef |
|---|---|---|---|---|---|---|---|---|---|---|---|
| **e** | e | ef | e.ef | efe | e.fef | fef | f | f.fe | fe | ef.fe | fe.ef |
| **f** | fe | f | fef | f.fe | efe | e.fef | fe.ef | ef | e.ef | e | ef.fe |
| **ef** | efe | fe | ef | ef.fe | fe.ef | f | e.ef | e | f.fe | e.fef | fef |
| **fe** | ef | fef | fe.ef | fe | e | ef.fe | e.fef | e.ef | f | f.fe | efe |
| **fef** | ef.fe | f.fe | fe | e.fef | fef | e.ef | ef | efe | e | fe.ef | f |
| **efe** | e.ef | fe.ef | e.fef | ef | f.fe | efe | ef.fe | f | fef | fe | e |
| **f.fe** | fe.ef | efe | e | fef | e.ef | fe | f.fe | e.fef | ef.fe | f | ef |
| **ef.fe** | e.fef | e | f.fe | fe.ef | f | ef | fe | ef.fe | efe | fef | e.ef |
| **fe.ef** | f | e.fef | ef.fe | e.ef | fe | e | efe | fef | fe.ef | ef | f.fe |
| **e.ef** | fef | ef.fe | efe | f | ef | f.fe | e | fe.ef | e.fef | e.ef | fe |
| **e.fef** | f.fe | e.ef | f | e | ef.fe | fe.ef | fef | fe | ef | efe | e.fef |

Let the ordered 19-tuple (1,2,3,4,5,6,7,8,9,10,11,12,13,14,15,16,17,18,19) =
(*e, f = e(ef.fe), ef, fe, fef, f.fe, e.ef, efe, e(f.fe) = f[e(e.ef)], ef.fe, fe.ef, f.efe, e.fef = f(ef.fe), efe.fef = f(f.fe), fef.efe = e(e.ef), e.efe = f.fef, (fe)(fe.ef), (ef)(ef.fe) = f(e.fef), f(e.ef)* ). Then the Cayley table for (5) *e = f(fe.ef)* is as follows:

**Table 9.** $e = f(fe.ef)$

|    | 1  | 2  | 3  | 4  | 5  | 6  | 7  | 8  | 9  | 10 | 11 | 12 | 13 | 14 | 15 | 16 | 17 | 18 | 19 |
|----|----|----|----|----|----|----|----|----|----|----|----|----|----|----|----|----|----|----|----|
| **1**  | 1  | 3  | 7  | 8  | 13 | 9  | 15 | 16 | 14 | 2  | 12 | 17 | 4  | 19 | 11 | 6  | 18 | 10 | 5  |
| **2**  | 4  | 2  | 5  | 6  | 16 | 14 | 19 | 12 | 8  | 13 | 1  | 3  | 18 | 10 | 9  | 7  | 11 | 17 | 15 |
| **3**  | 8  | 4  | 3  | 10 | 11 | 15 | 13 | 9  | 19 | 18 | 14 | 1  | 6  | 16 | 5  | 17 | 7  | 12 | 2  |
| **4**  | 3  | 5  | 11 | 4  | 19 | 12 | 14 | 10 | 1  | 15 | 17 | 7  | 2  | 8  | 16 | 18 | 13 | 6  | 9  |
| **5**  | 10 | 6  | 4  | 12 | 5  | 13 | 11 | 15 | 16 | 7  | 19 | 8  | 17 | 9  | 2  | 1  | 3  | 14 | 18 |
| **6**  | 11 | 16 | 19 | 5  | 1  | 6  | 17 | 4  | 7  | 12 | 18 | 13 | 9  | 3  | 8  | 15 | 14 | 2  | 10 |
| **7**  | 16 | 10 | 8  | 9  | 3  | 18 | 7  | 2  | 11 | 17 | 13 | 19 | 12 | 5  | 4  | 14 | 1  | 15 | 6  |
| **8**  | 7  | 11 | 13 | 3  | 14 | 10 | 12 | 8  | 17 | 9  | 6  | 18 | 5  | 1  | 19 | 2  | 15 | 4  | 16 |
| **9**  | 12 | 17 | 6  | 14 | 2  | 11 | 5  | 13 | 9  | 3  | 16 | 10 | 1  | 15 | 18 | 8  | 4  | 19 | 7  |
| **10** | 13 | 19 | 14 | 11 | 17 | 4  | 6  | 3  | 18 | 10 | 2  | 15 | 16 | 7  | 1  | 9  | 12 | 5  | 8  |
| **11** | 9  | 12 | 10 | 15 | 4  | 7  | 3  | 18 | 5  | 1  | 11 | 16 | 14 | 2  | 6  | 19 | 8  | 13 | 17 |
| **12** | 14 | 1  | 17 | 19 | 18 | 5  | 2  | 11 | 15 | 4  | 9  | 12 | 8  | 13 | 7  | 10 | 6  | 16 | 3  |
| **13** | 2  | 15 | 9  | 18 | 10 | 1  | 8  | 17 | 4  | 19 | 3  | 5  | 13 | 6  | 12 | 11 | 16 | 7  | 14 |
| **14** | 17 | 7  | 18 | 1  | 15 | 16 | 9  | 19 | 12 | 5  | 10 | 6  | 3  | 14 | 13 | 4  | 2  | 8  | 11 |
| **15** | 6  | 18 | 2  | 17 | 9  | 19 | 16 | 14 | 10 | 11 | 8  | 4  | 7  | 12 | 15 | 3  | 5  | 1  | 13 |
| **16** | 15 | 14 | 12 | 13 | 6  | 3  | 4  | 7  | 2  | 8  | 5  | 9  | 19 | 18 | 17 | 16 | 10 | 11 | 1  |
| **17** | 19 | 8  | 1  | 16 | 7  | 2  | 18 | 5  | 13 | 6  | 15 | 14 | 10 | 11 | 3  | 12 | 17 | 9  | 4  |
| **18** | 5  | 9  | 16 | 2  | 8  | 17 | 1  | 6  | 3  | 14 | 7  | 11 | 15 | 4  | 10 | 13 | 19 | 18 | 12 |
| **19** | 18 | 13 | 15 | 7  | 12 | 8  | 10 | 1  | 6  | 16 | 4  | 2  | 11 | 17 | 14 | 5  | 9  | 3  | 19 |

Let the ordered 19-tuple (1,2,3,4,5,6,7,8,9,10,11,12,13,14,15,16,17,18,19) =
= (e, f = e(fe.ef), ef, fe, fef, f.fe = e(e.ef), e.ef = f(f.fe) = f[e(e.ef)], efe = f(e.fef), e(f.fe), ef.fe,
fe.ef, f.efe, f.fef = (ef)(ef.fe), e.efe = (fe)(fe.ef), e.fef, f(e.ef) = e(ef.fe), e(f.fef), fef.efe, efe.fef ).
Then the Cayley table for (6) $e = f(ef.fe)$ is as follows:

**Table 10.** $e = f(ef.fe)$

|    | 1  | 2  | 3  | 4  | 5  | 6  | 7  | 8  | 9  | 10 | 11 | 12 | 13 | 14 | 15 | 16 | 17 | 18 | 19 |
|----|----|----|----|----|----|----|----|----|----|----|----|----|----|----|----|----|----|----|----|
| 1  | 1  | 3  | 7  | 8  | 15 | 9  | 6  | 14 | 4  | 16 | 2  | 5  | 17 | 11 | 13 | 18 | 19 | 12 | 10 |
| 2  | 4  | 2  | 5  | 6  | 13 | 7  | 16 | 12 | 19 | 1  | 9  | 14 | 10 | 17 | 8  | 3  | 18 | 11 | 15 |
| 3  | 8  | 4  | 3  | 10 | 11 | 18 | 15 | 9  | 6  | 13 | 19 | 2  | 14 | 5  | 17 | 1  | 16 | 7  | 12 |
| 4  | 3  | 5  | 11 | 4  | 16 | 12 | 19 | 10 | 2  | 18 | 14 | 17 | 8  | 13 | 1  | 7  | 9  | 15 | 6  |
| 5  | 10 | 6  | 4  | 12 | 5  | 1  | 11 | 18 | 15 | 17 | 16 | 19 | 9  | 2  | 14 | 8  | 13 | 3  | 7  |
| 6  | 11 | 13 | 16 | 5  | 9  | 6  | 14 | 4  | 17 | 12 | 8  | 1  | 3  | 18 | 7  | 15 | 10 | 19 | 2  |
| 7  | 14 | 10 | 8  | 9  | 3  | 13 | 7  | 16 | 12 | 5  | 15 | 6  | 19 | 4  | 2  | 17 | 11 | 1  | 18 |
| 8  | 7  | 11 | 15 | 3  | 19 | 10 | 2  | 8  | 5  | 9  | 17 | 13 | 1  | 16 | 18 | 12 | 14 | 6  | 4  |
| 9  | 2  | 14 | 17 | 19 | 1  | 11 | 18 | 15 | 9  | 3  | 12 | 10 | 6  | 8  | 4  | 5  | 7  | 13 | 16 |
| 10 | 15 | 16 | 19 | 11 | 14 | 4  | 17 | 3  | 13 | 10 | 1  | 18 | 7  | 9  | 12 | 6  | 8  | 2  | 5  |
| 11 | 9  | 12 | 10 | 18 | 4  | 17 | 3  | 13 | 7  | 2  | 11 | 15 | 16 | 6  | 19 | 14 | 5  | 8  | 1  |
| 12 | 19 | 9  | 14 | 16 | 8  | 5  | 1  | 11 | 18 | 4  | 7  | 12 | 15 | 10 | 6  | 2  | 3  | 17 | 13 |
| 13 | 18 | 7  | 12 | 1  | 6  | 14 | 4  | 17 | 3  | 19 | 5  | 11 | 13 | 15 | 16 | 9  | 2  | 10 | 8  |
| 14 | 6  | 19 | 2  | 15 | 17 | 3  | 13 | 7  | 16 | 8  | 18 | 9  | 12 | 14 | 10 | 4  | 1  | 5  | 11 |
| 15 | 16 | 18 | 9  | 13 | 10 | 2  | 8  | 5  | 1  | 6  | 3  | 7  | 11 | 12 | 15 | 19 | 4  | 14 | 17 |
| 16 | 13 | 1  | 18 | 17 | 12 | 19 | 10 | 2  | 8  | 15 | 4  | 3  | 5  | 7  | 11 | 16 | 6  | 9  | 14 |
| 17 | 12 | 15 | 6  | 7  | 2  | 8  | 5  | 1  | 11 | 14 | 13 | 16 | 18 | 19 | 9  | 10 | 17 | 4  | 3  |
| 18 | 17 | 8  | 1  | 14 | 7  | 16 | 12 | 19 | 10 | 11 | 6  | 4  | 2  | 3  | 5  | 13 | 15 | 18 | 9  |
| 19 | 5  | 17 | 13 | 2  | 18 | 15 | 9  | 6  | 14 | 7  | 10 | 8  | 4  | 1  | 3  | 11 | 12 | 16 | 19 |

Note that the one-step quasigroups determined by the values of the element $e$ in (1) through (6) satisfy the following identities:

(1) $e = f[f(fef)]$ satisfies the identity $x = y[y(yxy)]$
(2) $e = f[f(ef \cdot f)]$ satisfies the identity $x = y[y(xy.y)]$
(3) $e = (f \cdot fef)(ef) = (fef)(fe)$ satisfies the identity $x = (yxy)(yx)$
(4) $e = (fef \cdot f)(ef) = fe.fef$ satisfies the identity $x = (yx)(yxy)$
(5) $e = f[(ef)(fef)] = f(fe.ef)$ satisfies the identity $x = y(yx.xy)$
(6) $e = f[(ef)(ef.f)] = f(ef.fe)$ satisfies the identity $x = y(xy.yx)$

Therefore, the quasigroups determined by the Cayley tables 9 and 10 are NOT isomorphic. If they were isomorphic then the quasigroups determined by Tables 9 and 10 would satisfy the identities $x = y(yx.xy) = y(xy.yx)$. By cancellation, $xy.yx = yx.xy$. This contradicts Proposition 4 when $x = e$ and $y = f$.

# 6. One-step *W*-quasigroups of dimension greater than four

The methods we have used thus far to determine all one-step *W*-quasigroups of dimension less than or equal to four may be able to be used to determine one-step *W*-quasigroups of dimension greater than four. The steps are as follows:

1. Determine all possible values of *e* when $\|e\| = n$ (*n* greater than 4).
2. Use Proposition A and that G is cancellative to calculate the Cayley table of G.
3. Check that G is right modular.
4. Check that *e* and *f* can be written as words in any two distinct elements of G.

**Definition.** We define M(*n*) to be the maximum number of elements of a one-step *W*-quasigroup of dimension *n*. (So we can see that M(2) = 4, M(3) = 11 and M(4) = 29.)

Even for groups, rings and semigroups, the structures for which all finite one-step non-commutative members are known, there are questions about infinite one-step stuctures. As mentioned in the introduction, there are infinite, one-step non-commutative semigroups. As far as we are aware Redei's question as to whether there is an infinite one-step non-commutative group remains open.

*Open questions:*   (1) *Is M(n) < M(n+1), for every positive integer n?*
(2) *Is every one-step W-quasigroup finite?*
(3) *Is there a computer programme that will determine all possible values of e when $\|e\| = n$ (n greater than 4)?*
(4) *Suppose that we know the Cayley table of an idempotent groupoid G of order n. Suppose that we randomly check k number of products xy.z (k less than n cubed) and xy.z = zy.x for all of them. What is the probability that G is right modular?*

**References.**

[1]10 Albert Mansions, Crouch Hill, London N8 9RE; bobmonzo@talktalk.net